\documentclass[11pt,a4paper]{article}

\usepackage{amsmath,amssymb}
\newcommand{\C}{\mathbb{C}}
\newcommand{\R}{\mathbb{R}}
\newcommand{\N}{\mathbb{N}}
\newcommand{\Z}{\mathbb{Z}}
\def\cA{{\mathcal A}}
\def\cC{{\mathcal C}}
\def\cD{{\mathcal D}}
\def\cS{{\mathcal S}}
\newcommand{\ee}{\varepsilon}
\renewcommand{\div}{{\rm div}\,}
\newcommand{\Supp}{{\rm Supp}\,}
\newcommand{\Int}{\displaystyle \int}
\newcommand{\Frac}{\displaystyle \frac}
\newcommand{\Inf}{\displaystyle \inf}
\newcommand{\Sup}{\displaystyle \sup}
\newcommand{\Lim}{\displaystyle \lim}
\newcommand{\Max}{\displaystyle \max}
\newcommand{\Min}{\displaystyle \min}
\newcommand{\Sum}{\displaystyle \sum}

\def\d{\partial}
\def\ddl{\dot \Delta_l}
\def\ddq{\dot \Delta_q}
\def\tilde{\widetilde}
\def\hat{\widehat}
\newcommand{\D}{\Delta}
\newcommand{\La}{\Lambda}
\newcommand{\n}{\nabla}
\newcommand{\fd}{\frac{d}{2}}
\newcommand{\e}{\epsilon}
\newcommand{\fdp}{\frac{d}{p}}
\newcommand{\p}{\partial}
\newcommand{\ql}{q_{l}}
\newcommand{\ul}{u_{l}}
\newcommand{\hl}{h_{l}}
\newcommand{\Fl}{F_{l}}
\newcommand{\Gl}{G_{l}}
\newcommand{\Rl}{R_{l}}
\newcommand{\Rlp}{R_{l}'}
\newcommand{\g}{\int_{\mathbb{R}^{d}}}

\newcommand{\h}{\hookrightarrow}
\newcommand{\w}{\widetilde{T}}
\newcommand{\de}{\delta}
\newcommand{\del}{\bar{\delta}}
\newcommand{\ka}{\bar{\kappa}}
\newcommand{\cd}{\overline}
\newcommand{\qe}{q_\ee}
\newcommand{\ue}{u_\ee}
\newcommand{\dq}{\delta q}
\newcommand{\du}{\delta u}

\newtheorem{theorem}{Theorem}
\newtheorem{remarka}{Remark}
\newtheorem{lem}{Lemma}
\newtheorem{definition}{Definition}

\newtheorem{proposition}{Proposition}

\newtheorem{rem}{Remark}

\def\thesection{\arabic{section}}
\def\theequation{\arabic{section}.\arabic{equation}}
\def\thetheorem{\arabic{section}.\arabic{theorem}}
\def\theproposition{\arabic{section}.\arabic{proposition}}
\def\thecorollary{\arabic{section}.\arabic{corollaire}}
\def\thedefinition{\arabic{section}.\arabic{definition}}
\def\theremark{\arabic{section}.\arabic{remarque}}
\newcommand{\reset}{\setcounter{equation}{equation}\setcounter{theorem}{theorem}
\setcounter{proposition}{0}\setcounter{corollary}{corollaire}
\setcounter{remark}{remark}}
\title{Existence of global strong solution and vanishing capillarity-viscosity limit in one dimension for the Korteweg system}

\author{Fr\'ed\'eric Charve\footnote{Universit\'e Paris-Est Cr\'eteil, Laboratoire d'Analyse et de Math\'ematiques Appliqu\'ees (UMR 8050), 61 Avenue du G\'en\'eral de Gaulle, 94 010 Cr\'eteil Cedex (France). E-mail: frederic.charve@univ-paris12.fr}, Boris Haspot \thanks{
Ceremade UMR CNRS 7534
Universit\'e de Paris IX- Dauphine,
Place du Mar\'echal DeLattre De Tassigny
75775 PARIS CEDEX 16 , haspot@ceremade.dauphine.fr} \footnote{Basque Center of Applied Mathematics, Bizkaia Technology Park, Building 500,
E-48160, Derio (Spain)}
}

\date{}
\begin{document}

\maketitle

\begin{abstract}
In the first part of this paper, we prove the existence of global strong solution for Korteweg system in one dimension. 
In the second part, motivated by the processes of  vanishing capillarity-viscosity limit in order to select the physically relevant solutions for a hyperbolic system, we show that the global strong solution of the Korteweg system converges in the case of a $\gamma$ law for the pressure ($P(\rho)=a\rho^{\gamma}$, $\gamma>1$) to entropic solution of the compressible Euler equations. In particular it justifies that the Korteweg system is suitable for selecting the physical solutions in the case where the Euler system is strictly hyperbolic. The problem remains open for a Van der Waals pressure because in this case the system is not strictly hyperbolic and in particular the classical theory of Lax and Glimm (see \cite{Lax,G})  can not be used.
%
\end{abstract}
\section{Introduction}
We are concerned with compressible fluids endowed with internal
capillarity. The model we consider  originates from the XIXth
century work by Van der Waals and Korteweg \cite{VW,fK} and was
actually derived in its modern form in the 1980s using the second
gradient theory, see for instance \cite{fDS,fJL,fTN}. The first investigations begin with the Young-Laplace theory which claims that the phases are separated by a hypersurface and that the jump in the pressure across the hypersurface is proportional to the curvature of the hypersurface. The main difficulty consists in describing the location and the movement of the interfaces.\\
Another major problem is to understand whether the interface behaves as a discontinuity in the state space (sharp interface) or whether the phase boundary corresponds to a more regular transition (diffuse interface, DI).
The diffuse interface models have the advantage to consider only one set of equations in a single spatial domain (the density takes into account the different phases) which considerably simplifies the mathematical and numerical study (indeed in the case of sharp interfaces, we have to treat a problem with free boundary).\\
Let us consider a fluid of density $\rho\geq 0$, velocity field $u\in\R$, we are now interested in the following
compressible capillary fluid model, which can be derived from a Cahn-Hilliard like free energy (see the
pioneering work by J.- E. Dunn and J. Serrin in \cite{fDS} and also in
\cite{fA,fC,fGP,HM}).
The conservation of mass and of momentum write:
\begin{equation}
\begin{cases}
\begin{aligned}
&\frac{\p}{\p t}\rho^{\e}+\p_{x}(\rho^{\e} u^{\e})=0,\\
&\frac{\p}{\p t}(\rho^{\e} u^{\e})+\p_{x}(\rho^{\e}
(u^{\e})^{2})-\e\p_{x}(\rho^{\e}\p_{x}u^{\e})+\p_{x}(a(\rho^{\e})^{\gamma})=\e^{2}\p_{x}K,
\end{aligned}
\end{cases}
\label{3systeme}
\end{equation}
where the Korteweg tensor reads as following:
\begin{equation}
{\rm div}K
=\p_{x}\big(\rho^{\e}\kappa(\rho^{\e})\p_{xx}\rho^{\e}+\frac{1}{2}(\kappa(\rho^{\e})+\rho^{\e}\kappa^{'}(\rho^{\e}))|\p_{x}\rho^{\e}|^{2}\big)
-\p_{x}\big(\kappa(\rho^{\e})(\p_{x}\rho^{\e})^{2}\big).
\label{divK}
\end{equation}
$\kappa$ is the coefficient of capillarity and is a regular function of the form $\kappa(\rho)=\e^{2}\rho^{\alpha}$ with $\alpha\in\R$. In the sequel we shall assume that $\kappa(\rho)=\frac{\e^{2}}{\rho}$. The term
$\p_{x}K$  allows to describe the variation of density at the interfaces between two phases, generally a mixture liquid-vapor. $P=a\rho^{\gamma}$ with $\gamma\geq 1$ is a general $\gamma$ law pressure term.
$\e$ corresponds to the controlling parameter on the amplitude of the viscosity and of the capillarity.
When we set $v^{\e}=u^{\e}+\e\p_{x}(\ln\rho^{\e})$, we can write (\ref{3systeme}) on the following form (we refer to \cite{Hprepa} for the computations):
\begin{equation}
\begin{cases}
\begin{aligned}
&\frac{\p}{\p t}\rho^{\e}+\p_{x}(\rho^{\e} v^{\e})-\e\p_{xx}\rho^{\e}=0,\\
&\frac{\p}{\p t}(\rho^{\e} v^{\e})+\p_{x}(\rho^{\e}
(u^{\e})(v^{\e}))-\e\p_{x}(\rho^{\e}\p_{x}v^{\e})+\p_{x}(a(\rho^{\e})^{\gamma})=0,
\end{aligned}
\end{cases}
\label{1.1}
\end{equation}
We now consider the Cauchy problem of (\ref{1.1}) when the fluid is away from vacuum. Namely, we shall study (\ref{1.1}) with the following initial data:
\begin{equation}
\rho^{\e}(0,x)=\rho^{\e}_{0}(x)>0,\;\;u^{\e}(0,x)=u^{\e}_{0}(x),
\label{1.2}
\end{equation}
such that:
$$\lim_{x\rightarrow +,-\infty}(\rho_{0}^{\e}(x),u_{0}^{\e}(x))=(\rho^{+,-},u^{+,-}),\;\mbox{with}\;\;\rho^{+,-}>0.$$
We would like to study in the sequel the limit process of system (\ref{1.1}) when $\e$ goes to $0$ and to prove in particular that we obtain entropic solutions of the Euler system:
\begin{equation}
\begin{cases}
\begin{aligned}
&\frac{\p}{\p t}\rho+\p_{x}(\rho v)=0,\\
&\frac{\p}{\p t}(\rho v)+\p_{x}(\rho v^{2})+\p_{x}(a \rho^{\gamma})=0,
\end{aligned}
\end{cases}
\label{1.3}
\end{equation}
Let us now explain the interest of the capillary solutions for the hyperbolic systems of conservation laws.
\subsection{Viscosity capillarity processes of selection for the Euler system}
In addition of modeling a liquid-vapour mixture, the Korteweg also shows purely theoretical interests consisting in the selection of the physically relevant solutions of the Euler model (in particular when the system is not strictly hyperbolic).
The typical case corresponds to a Van der Waals pressure: indeed in this case the system is not strictly hyperbolic in the elliptic region (which corresponds to the region where the phase change occurs).
\\
In the adiabatic pressure framework ($P(\rho)=\rho^\gamma$ with $\gamma>1$), the system  is strictly hyperbolic and the theory is classical. More precisely we are able to solve the Riemann problem when the initial Heaviside data is small in the BV space. Indeed we are in the context of the well known Lax result as the system is also genuinely nonlinear (we refer to \cite{Lax}). It means we have existence of global $\mathcal{C}^1$-piecewise solutions which are unique in the class of the entropic solutions.
\\
This result as been extent by Glimm in the context of small initial data in the BV-space by using a numerical scheme and approximating the initial BV data by a $\mathcal{C}^1$-piecewise function (which implies to locally solve the Riemann problem via the Lax result). For the uniqueness of the solution we refer to the work of Bianchini and Bressan (\cite{BB1}) who use a viscosity method.
\\
In the setting of the Van der Waals pressure, the existence of global solutions and the nature of physical relevant solutions remain completely open. Indeed the system is not strictly hyperbolic anymore.


If we rewrite the compressible Euler system in Lagrangian coordinates by using the specific volume $\tau=1/\rho$ in $(\frac{1}{b},\infty)$ and the velocity $u$, the system satisfies in $(0,+\infty)\times\R$ the equations:
\begin{equation}
\begin{aligned}
&\p_{t}\tau-\p_{x}u=0,\\
&\p_{t}u-\p_{x}(\widetilde{P}(\tau))=0,
\end{aligned}
\label{euler}
\end{equation}
with the function $\widetilde{P}:(\frac{1}{b},\infty)\rightarrow (0,\infty)$ given by:
$$\widetilde{P}(\tau)=P(\frac{1}{\tau}),\;\;\;\tau\in(\frac{1}{b},\infty).$$

The two eigenvalues of the system are:
\begin{equation}
\lambda_{1}(\tau,v)=-\sqrt{-\widetilde{P}^{'}(\tau)},\;\;\;\lambda_{2}(\tau,v)=-\sqrt{-\widetilde{P}^{'}(\tau)}.
\label{vp}
\end{equation}
The corresponding eigenvectors $r_{1}$, $r_{2}$ are:
\begin{equation}
w_{1}(\tau,v)=\left(\begin{array}{c}
1\\
\sqrt{-\widetilde{P}^{'}(\tau)}\\
\end{array}
\right),\;\;w_{2}(\tau,v)=\left(\begin{array}{c}
1\\
-\sqrt{-\widetilde{P}^{'}(\tau)}\\
\end{array}
\right)
\end{equation}
Furthermore by calculus we obtain:
\begin{equation}
\n\lambda_{1}(\tau,v)\cdot w_{1}(\tau,v)=\frac{\widetilde{P}^{''}(\tau)}{2\sqrt{-\widetilde{P}^{'}(\tau)}},\;\;
\n\lambda_{2}(\tau,v)\cdot w_{2}(\tau,v)=\frac{-\widetilde{P}^{''}(\tau)}{2\sqrt{-\widetilde{P}^{'}(\tau)}}
\end{equation}
We now recall the definition of a \textit{standard conservation law} in the sense of Lax (it means entropy solutions):
\begin{itemize}
\item The system is \textbf{strictly entropic} if the eigenvalues are distinct and real.
\item The characteristics fields are \textbf{genuinely nonlinear} if we have for all $(\tau,v)$,
$$\n\lambda_{1}(\tau,v)\cdot w_{1}(\tau,v)\ne 0\;\;\mbox{and}\;\;\n\lambda_{2}(\tau,v)\cdot w_{2}(\tau,v)\ne 0,$$
\end{itemize}
for more details we refer to \cite{Serre}. The definition of genuine nonlinearity is some kind of extension of the notion of convexity to vector-valued functions (in particular when we consider the specific case of the traveling waves).  The previous assumptions aim at ensuring the existence and the uniqueness of the Riemann problem ( see \cite{Evans} and \cite{Serre}).

When $P$ is a Van der Waals pressure, we observe that the first conservation law (\cite{Serre,Evans}) is far from being a standard hyperbolic system, indeed:
\begin{itemize}
\item  It is not hyperbolic (but elliptic) in $(\frac{1}{\alpha_{1}},\frac{1}{\alpha_{2}})\times\R$,
\item the characteristic fields are not genuinely nonlinear in the hyperbolic part of the state space.
\end{itemize}
Here the classical Lax-Glimm theory cannot be applied. In particular there doesn't exist any entropy-flux pair, which suggests that the entropy framework is not adapted for selecting the physically relevant solutions. In order to deal with this problem, Van der Waals and Korteweg began by considering the stationary problem with null velocity, and solving $\nabla P(\rho)=0$. For more details we refer to \cite{Rohdehdr}. It consists in minimizing in the following admissible set
$$
A_{0}=\{\rho\in L^{1}(\Omega)/ W(\rho)\in L^{1}(\Omega), \;\int_{\Omega}\rho(x)dx=m\},
$$
the following functionnal
$$F[\rho]=\int_{\Omega}W(\rho(x))dx.$$
Unfortunately this minimization problem has an infinity of solutions, and many of them are physically irrelevant. In order to overcome this difficulty, Van der Waals in the XIX-th century was the first to regularize the previous functional by adding a quadratic term in the density gradient. More precisely he considered the following functional:
$$F^{\e}_{local}=\int_{\Omega}\big(W(\rho^{\e}(x))+\gamma\frac{\e^{2}}{2}|\n\rho^{\e}|^{2}\big)dx,$$
with:
$$A_{local}=H^{1}(\Omega)\cap A_{0}.$$
This variational problem has a unique solution and its limit (as $\e$ goes to zero) converge to a physical solution of the equilibrium problem for the Euler system with Van der Waals pressure, that was proved by Modica in \cite{REF} with the use of gamma-convergence.

By the Euler-Lagrange principle, the minimization of the Van der Waals functional consists in solving the following stationary problem:
$$
\nabla P(\rho^\e)=\gamma \e^2 \rho^\e \nabla \Delta \rho^\e,
$$
where the right-hand side can be expressed as the divergence of the capillarity tensor.

Heuristically, we also hope that the process of vanishing capillarity-viscosity limit selects the physical relevant solutions as it does for the stationary system. This problem actually remains open.

\subsection{Existence of global entropic solutions for Euler system}

Before presenting the results of this paper let us recall the results on this topic in these last decades. We shall focus on the case of a $\gamma$ pressure law $P(\rho)=a\rho^{\gamma}$ with $\gamma>1$ and $a$ positive. Let us mention that these cases are the only ones well-known (essentially because the system is strictly hyperbolic in this case and that we can exhibit many entropy-flux pairs). Here the Lax-Glimm theory can be applied, however at the end of the 70's, one was interested in relaxing the conditions on the initial data by only assuming $\rho_0$ and $u_0$ in $L^\infty$.

In the beginning of the 80's Di Perna initiated this program, consisting in obtaining global entropic solutions for $L^\infty$ initial data.


Indeed in \cite{Di1, Di2}, Di Perna prove the existence of global weak entropy solution of (\ref{1.3}) for $\gamma=1+\frac{2}{2d+1}$ and $\gamma=2k+\frac{3}{2k}+1$ (with $k\geq 1$), $d\geq 2$ by using the so-called "compensated compactness" introduced by Tartar in \cite{Ta}.
This result was extended by Chen in  \cite{Chen} in the case $\gamma\in(1,\frac{5}{3}]$ and by Lions et al in \cite{35} in the case $\gamma\in [3,\infty)$. In \cite{36}, Lions et al generalize this result to the general case $\gamma\in (1,3)$, and finally the case $\gamma=1$ is treated by \cite{Hu1}. We would like to mention that these results are obtained through a vanishing artificial viscosity on both density and velocity.\\
The problem of vanishing physical viscosity limit of compressible Navier-Stokes equations to compressible Euler equations was until recently an open problem. However Chen and Perepelista in \cite{10} proved that the solutions of the compressible Navier-Stokes system with constant viscosity coefficients converge to a entropic solution of the Euler system with finite energy. This result was extended in \cite{Hu2} to the case of viscosity coefficients depending on the density.\\
Inspired by \cite{10} and \cite{Hu2}, we would like to show that the solution of the Korteweg system (\ref{1.1}) converges to a entropic solution of the Euler system with finite energy when the pressure is a $\gamma$ law. To do this, we will prove for the first time up our knowledge the existence of global strong solution for the Korteweg system in one dimension in the case of Saint-Venant viscosity coefficients. By contrast, the problem of global strong solutions for compressible Navier-Stokes equations remains open (indeed one of the main difficulties consists in controlling the vacuum).
This result justifies that the Korteweg system allows us to select the relevant physical solutions of the compressible Euler system at least when the pressure is adiabatic ($P(\rho)=a \rho^\gamma$ with $\gamma>1$). The problem remains open in the case of a Van der Waals pressure.
\subsection{Results}
Let us now describe our main result. In the first theorem we prove the existence of global strong solution for the Korteweg system (\ref{1.1}).
\begin{theorem}
Let $\bar{\rho}>0$. Assume that the initial data $\rho_{0}$ and $u_{0}$ satisfy:
\begin{equation}
0<m_{0}\leq \rho_{0}\leq M_{0} <+\infty,  \rho_{0}-\bar{\rho}\in H^{1}(\R), v_{0}\in H^{1}(\R)\cap L^{\infty}(\R).
\label{2.5}
\end{equation}
Then there exists a global strong solution $(\rho, v)$ of (\ref{1.1}) on $\R^{+}\times\R$ such that for every $T>0$:
$$
\begin{aligned}
&\rho-\bar{\rho}\in L^{\infty}(0,T,H^{1}(\R)), \rho \in L^{\infty}(0,T,L^{\infty}(\R)), \\
&v\in L^{\infty}(0,T,H^{1}(\R))\cap L^{2}(0,T,H^{2}(\R))\;\;\mbox{and}\;\;v \in L^{\infty}(0,T,L^{\infty}(\R)).
\end{aligned}
$$
Finally this solution is unique in the class of weak solutions satisfying the usual energy inequality.
\label{theo}
\end{theorem}
\begin{rem}
We would like to point out that the problem remains open in the case of the Saint-Venant system, which corresponds to system (\ref{1.1}) without capillarity.
\end{rem}
In the following theorem, we are interested in proving the convergence of the global solutions of system (\ref{1.1}) to entropic solutions of the Euler system (\ref{1.3}).
\begin{theorem}
Let $\gamma>\frac{5}{3}$ and $(\rho^{\e},v^{\e})$ with $m_{\e}=\rho^{\e}v^{\e}$ be the global solution of the Cauchy problem (\ref{1.1}) with initial data $(\rho_{0}^{\e},v_{0}^{\e})$ as in theorem (\ref{theo}).
Then, when $\e\rightarrow 0$, there exists a subsequence of  $(\rho^{\e},m^{\e})$ that converge almost everywhere to a finite entropy solution $(\rho,\rho v)$ to the Cauchy problem (\ref{1.3}) with initial data $(\rho_{0},\rho_{0} v_{0})$.
\label{theo1}
\end{theorem}
\begin{remarka}
We would like to point out that Lions et al in \cite{36}  had obtained the existence of global entropic solution for $\gamma>1$ by a viscosity vanishing process, and the considered regularizing system was exactly the Korteweg system modulo the introduction of the effective velocity.
\end{remarka}
One important basis of our problem for theorem \ref{theo1} is the following compactness theorem established in \cite{10}.
\begin{theorem}(\textbf{Chen-Perepelitsa \cite{10}})
Let $\psi\in C^{2}_{0}(\R )$, $(\eta^{\psi},q^{\psi})$ be  a weak entropy pair generated by $\psi$. Assume that the sequences $(\rho^{\e}(x,t),v^{\e}(x,t))$ defined on $\R\times\R_{+}$ with $m^{\e}=\rho^{\e}v^{\e}$, satisfies the following conditions:
\begin{enumerate}
\item For any $-\infty<a<b<+\infty$ and all $t>0$, it holds that:
\begin{equation}
\int^{t}_{0}\int^{b}_{a}(\rho^{\e})^{\gamma+1}dxd\tau\leq C(t,a,b),
\label{1.8}
\end{equation}
where $C(t)>0$ is independent of $\e$.
\item For any compact set $K\subset\R$, it holds that
\begin{equation}
\int^{t}_{0}\int_{K}\big((\rho^{\e})^{\gamma+\theta}+\rho^{\e}|v^{\e}|^{3}\big)dxd\tau\leq C(t,K),
\label{1.9}
\end{equation}
where $C(t,K)>0$ is independent of $\e$.
\item The sequence of entropy dissipation measures
\begin{equation}
\eta^{\psi}(\rho^{\e},m^{\e})_{t}+q^{\psi}(\rho^{\e},m^{\e})_{x}\;\;\mbox{are compact in}\;\;H^{-1}_{loc}(\R^{2}_{+}).
\label{1.10}
\end{equation}
Then there is a subsequence of $(\rho^{\e},m^{\e})$ (still denoted $(\rho^{\e},m^{\e})$) and a pair of measurable functions $(\rho,m)$ such that:
\begin{equation}
(\rho^{\e},m^{\e})\rightarrow (\rho,m),\,\,\mbox{a.e as}\,\,\e\rightarrow 0.
\label{1.11}
\end{equation}
\end{enumerate}
\label{theo2}
\end{theorem}
\begin{remarka}
We would like to recall that the estimate (\ref{1.9}) was first derived by Lions et al in \cite{35} by relying the moment lemma introduced by Perthame in \cite{Per}.
\end{remarka}
The paper is arranged as follows. In section \ref{section2} we recall some important results on the notion of entropy enrtopy-flux pair for Euler system and on the kinetic formulation of Lions et al in \cite{35}. In section \ref{section3}, we show theorem \ref{theo} and in the last section \ref{section4} we prove theorem \ref{theo1}.
\section{Mathematical tools}
\label{section2}
\begin{definition}
 A pair of functions $(\eta(\rho,v), H(\rho,v))$ or $(\eta(\rho,m),q(\rho,m))$ for $m=\rho v$, is called an entropy-entropy flux pair of system (\ref{1.1}), if the following holds:
$$[\eta(\rho,v)]_{t}+[H(\rho,v)]_{x}=0,$$
for any smooth solution of (\ref{1.3}).
Furthermore $(\eta(\rho,v)$ is called a weak entropy if:
$$\eta(0,u)=0,\;\;\mbox{for any fixed}\;v.$$
\end{definition}
\begin{definition}
An entropy $\eta(\rho,m)$ is convex if the Hessian $\n^{2}\eta(\rho,m)$ is nonnegative definite in the region under consideration.
\end{definition}
Such $\eta$ satisfy the wave equation:
$$\p_{tt}\eta=\theta^{2}\rho^{\gamma-3}\p_{xx}\eta.$$ From \cite{35}, we obtain an explicit representation of any weak entropy $(\eta,q)$ under the following form:
\begin{equation}
\begin{aligned}
&\eta^{\psi}(\rho,m)=\int_{\R}\chi(\rho, s-v)\psi(s)ds,\\
&H^{\psi}(\rho,m)=\int_{\R}(\theta s+(1-\theta)u)\chi(\rho, s-v)\psi(s)ds,
\end{aligned}
\label{1.4}
\end{equation}
where the kernel $\chi$ is defined as follows:
$$\chi(\rho,v)=[\rho^{2\theta}-v^{2}]_{+}^{\lambda},\;\lambda=\frac{3-\gamma}{2(\gamma-1)}>-\frac{1}{2},\;\;\mbox{and}\;\theta=\frac{\gamma-1}{2},$$
and here:
$$
\begin{aligned}
t^{\lambda}_{+}=&t^{\lambda}\;\;\mbox{for}\;t>0,\\
=&0\;\;\mbox{for}\;t\leq0,
\end{aligned}
$$
\begin{proposition}
(see \cite{35})
\end{proposition}
For instance, when $\psi(s)=\frac{1}{2}s^{2}$, the entropy pair is the mechanical energy and the associated flux:
\begin{equation}
\eta^{*}(\rho,m)=\frac{m^{2}}{2\rho}+e(\rho),\;\;q^{*}(\rho,m)=\frac{m^{3}}{2\rho^{2}}+e^{'}(\rho),
\label{1.5}
\end{equation}
where $e(\rho)=\frac{\kappa}{\gamma-1}\rho^{\gamma}$ represents the gas internal energy in physics.\\
\\
In the sequel we will work far away of the vacuum that it why we shall introduce equilibrium states such that we avoid the vacuum. Let $(\bar{\rho}(x),\bar{v}(x))$ be a pair of smooth monotone functions satisfying $(\bar{\rho}(x),\bar{v}(x))=(\rho^{-,+},v^{-,+})$
when $-+x\geq L_{0}$ for some large $L_{0}>0$. The total mechanical energy for (\ref{1.1}) in $\R$ with respect to the pair of reference function $(\bar{\rho}(x),\bar{v}(x))$ is:
\begin{equation}
E[\rho,v](t)=\int_{\R}\big(\frac{1}{2}\rho(t,x)|v(t,x)-\bar{v}(x)|^{2}+e^{*}(\rho(t,x),\bar{\rho}(x))\big)dx
\label{1.7}
\end{equation}
where $e^{*}(\rho,\bar{\rho})=e(\rho)_e(\bar{\rho})-e^{'}(\bar{\rho}(\rho-\bar{\rho})\geq 0$.
The total mechanical energy for system (\ref{3system}) with $\kappa(\rho)=\frac{\kappa}{\rho}$ is:
\begin{equation}
E_{1}[\rho,u](t)=\int_{\R}\big(\frac{1}{2}\rho(t,x)|u(t,x)-\bar{u}(x)|^{2}+e^{*}(\rho(t,x),\bar{\rho}(x))+
\e^{2} (\p_{x}\rho^{\frac{1}{2}})^{2}\big)dx
\label{1.7}
\end{equation}
and the total mechanical energy for system (\ref{1.1}) is:
\begin{equation}
E_{2}[\rho,v](t)=\int_{\R}\big(\frac{1}{2}\rho(t,x)|v(t,x)-\bar{v}(x)|^{2}+e^{*}(\rho(t,x),\bar{\rho}(x))\big)dx
\label{1.7}
\end{equation}
\begin{definition}
Let $(\rho_{0},v_{0})$ be given initial data with finite-energy with respect to the end states: $( \rho^{\pm}, v^{\pm})$ at infinity, and $E[\rho_{0},v_{0}]\leq E_{0}<+\infty$. A pair of measurable functions $(\rho,u):\R^{2}_{+}\rightarrow \R^{2}_{+}$ is called a finite-energy entropy solution of the Cauchy problem (\ref{1.3}) if the following properties hold:
\begin{enumerate}
\item The total energy is bounded in time such that there exists a bounded function $C(E,t)$, defined on $\R^{+}\times\R^{+}$ and continuous in $t$ for each $E\in\R^{+}$ with for a.e $t>0$:
    $$E[\rho,v](t)\leq C(E_{0},t).$$
\item The entropy inequality:
$$\eta^{\psi}(\rho,v)_{t}+q^{\psi}(\rho,v)_{x}\leq 0,$$
is satisfied in the sense of distributions for all test functions $\psi(s)\in\{\pm 1,\pm s, s^{2}\}$.
\item The initial data $(\rho_{0},v_{0})$ are obtained in the sense of distributions.
\end{enumerate}
\end{definition}
We now give our main conditions on the initial data (\ref{1.2}), which is inspired from \cite{10}.
\begin{definition}
Let $(\bar{\rho}(x),\bar{v}(x))$ be some pair of smooth monotone functions satisfying $(\bar{\rho}(x),\bar{v}(x))=(\rho^{-,+},v^{-,+})$
when $-+x\geq L_{0}$ for some large $L_{0}>0$. For positive constant $C_{0}$, $C_{1}$ and $C_{2}$ independent of $\e$, we say that the initial data $(\rho_{0}^{\e},v_{0}^{\e})$ satisfy the condition ${\cal H}$ if they verify the following properties:
\begin{itemize}
\item $\rho_{0}^{\e}>0$, $\int_{\R}\rho_{0}^{\e}(x)|u_{0}^{\e}(x)-\bar{u}(x)|\leq C_{0}<+\infty,$
\item The energy is finite:
$$\int_{\R}\big(\frac{1}{2}\rho_{0}^{\e}(x)|v_{0}^{\e}(x)-\bar{v}(x)|^{2}+e^{*}(\rho_{0}^{\e}(x),\bar{\rho}(x))\big)dx\leq C_{1}<+\infty,
$$
\item $$\e^{2}\int_{\R}\frac{|\p_{x}\rho_{0}^{\e}(x)|^{2}}{\rho_{0}^{\e}(x)^{3-2\alpha}}dx\leq C_{2}<+\infty.$$
\end{itemize}
\end{definition}

In this section, we would like to recall some properties on the pair of entropy for the system (\ref{1.3}). Smooth solutions of (\ref{1.3}) satisfy the conservation laws:
$$\p_{t}\eta(\rho,u)+\p_{x}H(\rho,u)=0,$$
if and only if:
\begin{equation}
\eta_{\rho\rho}=\frac{P^{'}(\rho)}{\rho^{2}}\eta_{uu}.
\label{ondes}
\end{equation}
We supplement the equation \ref{ondes} by giving initial conditions:
\begin{equation}
\eta(0,u)=0,\\
\eta_{\rho}(0,u)=\psi(u).
\label{initial}
\end{equation}
We are now going to give a sequel of proposition on the properties of $\eta$, we refer to \cite{35} for more details.
\begin{proposition}
\label{proputile}
For $\rho\geq 0$, $u, \omega\in\R$,
\begin{itemize}
\item The fundamental solution of (\ref{ondes})-(\ref{initial}) is the solution corresponding to $\eta_{\rho}(0,u)=\delta(u)$ is given by:
\begin{equation}
\begin{aligned}
\chi(\rho,\omega)=(\rho^{\gamma-1}-\omega^{2})^{\lambda}_{+}\;\;\mbox{with}\;\;\lambda=\frac{3-\gamma}{2(\gamma-1)}.
\end{aligned}
\end{equation}
\item The solution of (\ref{ondes})-(\ref{initial}) is given by:
\begin{equation}
\eta(\rho,u)=\int_{\R}\psi(\xi)\chi(\rho,\xi-u)d\xi,
\label{entropie}
\end{equation}
\item $\eta$ is convex in $(\rho,\rho u)$ for all $\rho$, $u$ if and only if $g$ is convex.
\item The entropy flux $H$ associated with $\eta$ is given by:
\begin{equation}
H(\rho,u)=\int_{\R}\psi(\xi)[\theta\xi+(1-\theta)\xi]\chi(\rho,\xi-u)d\xi\;\;\mbox{where}\;\;\theta=\frac{\gamma-1}{2}.
\label{flux}
\end{equation}
\end{itemize}
\end{proposition}
We now give a important result on the entropy pair (see \cite{35}, lemma 4) .
\begin{proposition}
\label{pair35}
Taking $\psi(s)=\frac{1}{2}s|s|$, then there exists a positive constant $C>0$, depending only on $\gamma>1$, such that the entropy pair $(\eta^{\psi},H^{\psi})$ satisfies:
\begin{equation}
\begin{aligned}
&|\eta^{\psi}(\rho,u)|\leq (\rho|u|^{2}+\rho^{\gamma}),\\
&H^{\psi}(\rho,u)\geq C^{-1}(\rho|u|^{3}+\rho^{\gamma+\theta}),\;\;\mbox{for all}\;\rho\geq 0\;\mbox{and}\;u\in\R,\\
&|\eta^{\psi}_{m}(\rho,u)|\leq (\rho|u|+\rho^{\theta}),\\
&|\eta^{\psi}_{m m}(\rho,u)|\leq C\rho^{-1}.
\end{aligned}
\label{2.37}
\end{equation}
\end{proposition}
We are now going to give recent results on the entropy pair $(\eta^{\psi},q^{\psi})$ generated by $\psi\in C^{2}_{0}(\R)$ (we refer to  \cite{10} for more details).
\begin{proposition}
\label{propChen}
For a $C^{2}$ function $\psi:\R\rightarrow\R$, compactly supported on the interval $[a,b]$, we have:
\begin{equation}
\mbox{supp}(\eta^{\psi}),\mbox{supp}(q^{\psi})\subset\{(\rho,m)=(\rho,\rho u):\;u+\rho^{\theta}\geq a, u-\rho^{\theta}\leq b\}:
\label{3.2}
\end{equation}
Furthermore, there exists a constant $C_{\psi}$ such that, for any $\rho\geq 0$ and $u\in\R$, we have:
\begin{itemize}
\item For $\gamma\in(1,3]$,
\begin{equation}
|\eta^{\psi}(\rho,m)|+|q^{\psi}(\rho,m)|\leq C_{\psi}\rho.
\label{3.3}
\end{equation}
\item For $\gamma\in(3,+\infty)$,
\begin{equation}
|\eta^{\psi}(\rho,m)|\leq C_{\psi}\rho,\;\;|q^{\psi}(\rho,m)|\leq C_{\psi}(\rho+\rho^{\theta+1}).
\label{3.4}
\end{equation}
\item If $\eta^{\psi}$ is considered as a function of $(\rho,m)$, $m=\rho u$ then
\begin{equation}
|\eta_{m}^{\psi}(\rho,m)|+|\rho \eta_{mm}^{\psi}(\rho,m)|\leq C_{\psi},
\label{3.5}
\end{equation}
and, if $\eta^{\psi}_{m}$ is considered as a function of $(\rho,u)$, then
\begin{equation}
|\eta_{m}^{\psi}(m,u)|+|\rho^{1-\theta} \eta_{m\rho}^{\psi}(\rho,\rho u)|\leq C_{\psi}.
\label{3.6}
\end{equation}
\end{itemize}
\end{proposition}
We now would like to express the kinetic formulation of (\ref{1.3}) introduced in (\cite{35}).
\begin{theorem}
Let $(\rho,\rho v)\in L^{\infty}(\R^{+},L^{1}(\R))$ have finite energy and $\rho\geq 0$, then it is an entropy solution of (\ref{1.3}) if and only if there exists a non-positive bounded measure $m$ on $\R^{+}\times\R^{2}$ such that the function $\chi(\rho,\xi-u)$ satisfies:
\begin{equation}
\p_{t}\chi+\p_{x}[(\theta\xi+(1-\theta)u)\chi]=\p_{\xi\xi}m(t,x,\xi).
\end{equation}
\label{cinetique}
\end{theorem}
\section{Proof of theorem \ref{theo}}
\label{section3}
We would like to start with recalling an important result due to Solonnikov (see \cite{Sol}). Let $\rho_{0}$ the initial density such that:
\begin{equation}
0<m_{0}\leq \rho_{0}\leq M_{0}<+\infty.
\label{initiald}
\end{equation}
When the viscosity coefficient $\mu(\rho)$ satisfies:
\begin{equation}
\mu(\rho)\geq c>0\;\;\mbox{for all} \rho\geq 0,
\label{visco}
\end{equation}
we have the existence of strong solution for small time. More exactly, we have:
\begin{proposition}
Let $(\rho_{0},v_{0})$ satisfy (\ref{initiald}) and assume that $\mu$ satisfies (\ref{visco}), then there exists $T_{0}>0$ depending on $m_{0}$, $M_{0}$, $\|\rho_{0}-\bar{\rho}\|_{H^{1}}$ and $\|v_{0}\|_{H^{1}}$ such that (\ref{1.1}) has a unique solution $(\rho,v)$ on $(0,T_{0})$ satisfying:
$$
\begin{aligned}
&\rho-\bar{\rho}\in L^{\infty}(H^{1}(\R),\;\;\p_{t}\rho\in L^{2}((0,T_{1})\times\R ),\\
&v\in L^{2}(0,T_{1},H^{2}(\R)),\;\;\p_{t}v\in L^{2}((0,T_{1})\times\R)
\end{aligned}
$$
for all $T_{1}<T_{0}$.
\end{proposition}
\begin{remarka}
The main point in this theorem is that the time of existence $T_{0}$ depends only of the norms of $\rho_{0}$ which gives us a low bounds on $T_{0}$ of the system (\ref{1.1}).
\end{remarka}
In view of this proposition, we see that if we introduce a truncated viscosity coefficient $\mu_{n}(\rho)$:
$$\mu_{n}(\rho)=\max(\rho,\frac{1}{n}),$$
then there exists approximated solutions $(\rho_{n},v_{n})$ defined for small time $(0,T_{0})$ of the system (\ref{1.1}). In order to prove theorem \ref{theo}
, we only have to show that $(\rho_{n},v_{n})$ satisfies the following bounds uniformly with respect to $n$ and $T$ large:
\begin{equation}
\begin{aligned}
&0<m_{0}\leq \rho_{n}\leq M_{0} <+\infty,  \;\;\forall t\in[0,T],\\
&\rho_{n}-\bar{\rho}\in L^{\infty}_{T}(H^{1}(\R)), \\
&v_{n}\in  L^{\infty}_{T}(H^{1}(\R)).
\end{aligned}
\end{equation}
We are going to follow the method of Lions et al in \cite{36}, indeed the main point is to prove that we can extend the notion of Riemann invariant or more precisely the kinetic formulation of proposition \ref{cinetique} to the system (\ref{1.1}). We recall that system  (\ref{1.1}) has the following form:
\begin{equation}
\begin{cases}
\begin{aligned}
&\frac{\p}{\p t}\rho_{n}+\p_{x}(\rho_{n} v_{n})-\e\p_{xx}\rho_{n}=0,\\
&\frac{\p}{\p t}(\rho_{n}v_{n})+\p_{x}(\rho_{n}
v_{n}v_{n})-\e\p_{x}(\p_{x}\rho_{n}v_{n})-\e\p_{x}(\rho_{n}\p_{x}v_{n})+\p_{x}(a(\rho_{n})^{\gamma}=0,
\end{aligned}
\end{cases}
\label{1.1a}
\end{equation}
and we have finally:
\begin{equation}
\begin{cases}
\begin{aligned}
&\frac{\p}{\p t}\rho_{n}+\p_{x}(\rho_{n} v_{n})-\e\p_{xx}\rho_{n}=0,\\
&\frac{\p}{\p t}(\rho_{n}v_{n})+\p_{x}(\rho_{n}
v_{n}v_{n})-\e\p_{x}\p_{x}(\rho_{n} v_{n})+\p_{x}(a(\rho_{n})^{\gamma}=0,
\end{aligned}
\end{cases}
\label{1.1b}
\end{equation}
Following \cite{36} and setting $m_{n}=\rho_{n}v_{n}$ we have for any pair of entropy flux $(\eta(\rho,u), H(\rho,u))$ defined by (\ref{entropie}) and (\ref{flux}) where $\eta$ is a convex function of $(\rho_{n},m_{n})$. We write $\eta=\bar{\eta}(\rho_{n},m_{n})$:
$$
\begin{aligned}
\p_{t}\eta
+\p_{x}H&=\e\bar{\eta}_{\rho}\p_{xx}\rho_{n}+\e\bar{\eta}_{m}\p_{xx}m_{n},\\
&=\e\p_{xx}\eta-\e(\bar{\eta}_{\rho\rho}(\p_{x}\rho_{n})^{2}+2\bar{\eta}_{\rho m}(\p_{x}\rho_{n})(\p_{x}m_{n})+\bar{\eta}_{mm}(\p_{x}m_{n})^{2}).
\end{aligned}
$$
Here we define $\mu_{n}$ such that:
$$\mu_{n}=\bar{\eta}_{\rho\rho}(\p_{x}\rho_{n})^{2}+2\bar{\eta}_{\rho m}(\p_{x}\rho_{n})(\p_{x}m_{n})+\bar{\eta}_{mm}(\p_{x}m_{n})^{2}$$
By proposition \ref{proputile}, we can check that $\mu_{n}\geq 0$. We obtain then that:
$$
\p_{t}\eta(\rho_{n},v_{n})+\p_{x}H((\rho_{n},v_{n})-\e\bar{\eta}_{\rho}\p_{xx}\rho_{n}\leq 0\;\;\mbox{in}\;\R\times (0,+\infty).
$$
By applying the same method than for proving the theorem \ref{cinetique}, we obtain the following kinetic formulation:
\begin{equation}
\p_{t}\chi+\p_{x}([\theta\xi+(1-\theta)v_{n}]\chi)-\p_{xx}\chi=\p_{\xi\xi}\bar{m}_{n}\;\;\mbox{on}\;\R^{2}\times(0,+\infty),
\label{riemann}
\end{equation}
where $\bar{m}_{n}$ is a nonpositive bounded measure on $\R^{2}\times(0,+\infty)$.  Finally we recover the classical maximum principle by multiplying (\ref{riemann}) by the convex functions $g(\xi)=(\xi-\xi_{0})_{+}$ and $g(\xi)=(\xi-\xi_{0})_{-}$ and integrating over $\R^{2}\times(0,+\infty)$. Indeed as we have that:
$$-C\leq\min_{x}(v_{0}-\rho_{0}^{\theta})\leq\max_{x}(v_{0}+\rho_{0}^{\theta})\leq C,$$
and that:
$$\mbox{supp}\xi=[v-\rho^{\theta},v+\rho[{\theta}].$$
For $\xi_{0}$ large enough, we can show that:
$$\mbox{supp}\xi_{0}\cap\mbox{supp}\chi=\emptyset.$$
We have obtain then that:
$$-C\leq\min_{x}(v_{0}-\rho_{0}^{\theta})\leq v_{n}-\rho_{n}^{\theta}\leq v_{n}+\rho_{n}^{\theta}\leq\max_{x}(v_{0}+\rho_{0}^{\theta})\leq C.$$
In particular we obtained that $\rho_{n}$ and $v_{n}$ are uniformly bounded in $L^{\infty}(0,T_{n}, L^{\infty}(\R))$ or:
\begin{equation}
\sup_{x\in\R,t\in(0,T_{n})}\big(|\rho_{n}(t,x)|+|v_{n}(t,x)|)\leq C_{0},
\label{imp2}
\end{equation}
\section{Proof of theorem \ref{theo1}}
\subsection{Uniform estimates for the solutions of (\ref{1.1})}
\label{section4}
First we assume that $(\rho^{\e},v^{\e})$ is the global solutions of Korteweg's equations (\ref{1.1}) constructed in theorem \ref{theo} and satisfying:
\begin{equation}
\rho^{\e}(t,x)\geq c^{\e}(t),\,\mbox{for some}\;c^{\e}(t)>0,
\label{2.1}
\end{equation}
and
\begin{equation}
\lim_{x\rightarrow \pm \infty}(\rho^{\e},v^{\e})(x,t)=(\rho^{\pm},u^{\pm}).
\label{2.2}
\end{equation}
Here we are working around a non constant state $(\bar{\rho},\bar{v})$ with:
$$
\lim_{x\rightarrow \pm \infty}(\bar{\rho},\bar{v})(x,t)=(\rho^{\pm},u^{\pm}).
$$
It is a simple extension of theorem \ref{theo}. Our goal is now to check the properties (\ref{1.8}),  (\ref{1.9}) and (\ref{1.10}) in order to use the theorem \ref{theo2} of Chen and Perepelista (see \cite{10}) in order to prove the theorem \ref{theo1}.\\
For simplicity, throughout this section, we denote $(\rho,v)=(\rho^{\e},v^{\e})$ and $C>0$ denote the constant independent of $\e$.\\
We start with recalling the inequality energy for system (\ref{1.1}), indeed by the introduction of the effective velocity we obtain new entropies (see \cite{Hprepa}).
\begin{lem}
Suppose that 
$E_{1}[\rho_{0},u_{0}]\leq E_{0}<+\infty$ for some $E_{0}>0$ independent of $\e$. It holds that:
\begin{equation}
\sup_{0\leq\tau\leq t} E_{1}[\rho, u](\tau)+\e\int^{t}_{0}\int_{\R}\rho u_{x}^{2} dxd\tau\leq C(t),
\label{2.3a}
\end{equation}
and:
\begin{equation}
\sup_{0\leq\tau\leq t}  E_{2}[\rho, v](\tau)+\e\int^{t}_{0}\int_{\R}\rho v_{x}^{2} dxd\tau+\e\int^{t}_{0}\int_{\R}\rho^{\gamma-2}\rho_{x}^{2} dxd\tau\leq C(t),
\label{2.3}
\end{equation}
where $C(t)$ depends on $E_{0}$, $t$, $\bar{\rho}$, and $\bar{u}$ but not on $\e$.
\label{lemma1}
\end{lem}
{\bf Proof:} It suffices to writes the energy inequalities for system (\ref{1.1}) and from system (\ref{1.3}). More exactly we have:
$$\frac{d}{dt}\int_{\R}(\eta(\rho,m)-\eta(\bar{\rho},\bar{\rho}\bar{u})dx+\e\int_{\R}\rho u_{x}^{2}dx=q(\rho^{-},m^{-})-q(\rho^{+},m^{+}),$$
with the entropy pair:
$$\eta(\rho,m)=\frac{m^{2}}{2\rho}+e(\rho),\;\;q(\rho,m)=\frac{m^{3}}{2\rho^{2}}+m e^{'}(\rho),$$
with $e(\rho)=\frac{a}{\gamma-1}\rho^{\gamma}$.
Since we have:
$$\e(\rho,\bar{\rho})\geq \rho(\rho^{\theta}-\bar{\rho}^{\theta})^{2},\;\theta=\frac{\gamma-1}{2},$$
we can classically bootstrap on the left hand-side the term $q(\rho^{-},m^{-})-q(\rho^{+},m^{+})$.
\begin{remarka}
Since vacuum could occur in our solution, the inequality
$$\int^{t}_{0}\int_{\R}\rho u_{x}^{2} dxd\tau\leq C(t),$$
in (\ref{2.3}) is much weaker than the corresponding one in \cite{10}. That is why lemma \ref{lemme2} will be more tricky to obtain.
\end{remarka}
The following higher order integrability estimate is crucial in compactness argument.
\begin{lem}
\label{lemme2}
If the conditions of lemma \ref{lemma1} hold, then for any $-\infty<a<b<+\infty$ and all $t>0$, it holds that:
\begin{equation}
\int^{t}_{0}\int^{b}_{a}\rho^{\gamma+1}dxd\tau\leq C(t,a,b),
\label{2.21}
\end{equation}
where $C(t)>0$ depends on $E_{0}$, $a$, $b$, $\gamma$, $t$, $\bar{\rho}$, $\bar{u}$ but not on $\e$.
\end{lem}
\begin{remarka}
The proof follows the same ideas than in the case of compressible Navier-Stokes equations when we wish to obtain a gain of integrability on the density. We refer to \cite{fL2} for more details. The proof is also inspired from Huang et al in \cite{Hu2}.
\end{remarka}
{\bf Proof.} Choose $\omega\in C^{\infty}_{0}(\R)$ such that:
$$0\leq\omega(x)\leq 1,\;\omega(x)=1\;\;\mbox{for}\;x\in[a,b],\;\mbox{and}\;{\rm supp}{\omega}=(a-1,b+1).$$
By the momentum equation of (\ref{1.1}) and by localizing, we have
\begin{equation}
(P(\rho)\omega)_{x}=-(\rho  u v\omega)_{x}+(P(\rho)+\rho u v)\omega_{x}-(\rho v)_{t}\omega+\e(\rho v_{x}\omega)_{x}-\e\rho v_{x}\omega_{x}.
\label{2.22}
\end{equation}
Integrating (\ref{2.22}) with respect to spatial variable over $(-\infty,x)$, we obtain:
\begin{equation}
P(\rho)\omega=-\rho u v\omega+\e(\rho v_{x}\omega)_{x}-(\int^{x}_{-\infty}\rho v\,\omega dy)_{t}+\int^{x}_{-\infty}[(\rho u v+P(\rho))\omega_{x}-\e\rho v_{x}\omega_{x}.
\label{2.23}
\end{equation}
Multiplying (\ref{2.23}) by $\rho\omega$, we have
\begin{equation}
\begin{aligned}
\rho P(\rho)\omega^{2}=&-\rho^{2}u v \omega^{2}+\e\rho^{2}v_{x}\omega^{2}-(\rho\omega\int^{x}_{-\infty}\rho v\,\omega dy)_{t}\\
&-(\rho u)_{x}\omega(\int^{x}_{-\infty}\rho u\,\omega dy)+
\rho\omega\int^{x}_{-\infty}[(\rho u v+P(\rho))\omega_{x}-\e\rho v_{x}\omega_{x}]dx,\\
=&\e\rho^{2}v_{x}\omega^{2}-(\rho\omega\int^{x}_{-\infty}\rho v\,\omega dy)_{t}-(\rho u\omega\int^{x}_{-\infty}\rho v\,\omega dy)_{x}\\
&+\rho u\omega_{x}\int^{x}_{-\infty}\rho v\,\omega dy+\rho\omega\int^{x}_{-\infty}[(\rho u v+P(\rho))\omega_{x}-\e\rho v_{x}\omega_{x}]dx,
\end{aligned}
\label{2.24}
\end{equation}
We now integrate (\ref{2.24}) over $(0,t)\times\R$ and we get:
\begin{equation}
\begin{aligned}
&\int^{t}_{0}\int_{\R}a\rho^{\gamma+1}\omega^{2}dxd\tau=\e\int^{t}_{0}\int_{\R}\rho^{2}v_{x}\omega^{2}-\int_{\R}(\rho\omega\int^{x}_{-\infty}\rho v\,\omega dy)dx\\
&+\int_{\R}(\rho_{0}\omega\int^{x}_{-\infty}\rho_{0} v_{0}\,\omega dy)dx+\int^{t}_{0}\int_{\R}\big(\rho u\omega_{x}\int^{x}_{-\infty}\rho v\,\omega dy\big)dx d\tau\\
&\hspace{2cm}+\int^{t}_{0}\int_{\R}\big(\rho\omega\int^{x}_{-\infty}[(\rho u v+P(\rho))\omega_{x}-\e\rho v_{x}\omega_{x}]dx\big)dx d\tau.
\end{aligned}
\label{2.25}
\end{equation}
Let
\begin{equation}
A=\{ x:\rho(t,x)\geq\underline{\rho}\},\;\;\mbox{where}\;\underline{\rho}=2\max(\rho+,\rho-),
\label{2.26}
\end{equation}
then we have the following estimates by (\ref{2.3}):
\begin{equation}
|A|\leq \frac{C(t)}{e^{*}(2\underline{\rho},\bar{\rho})}=d(t).
\label{2.27}
\end{equation}
By (\ref{2.26}), for any $(t,x)$ there exists a point $x_{0}=x_{0}(t,x)$ such that $|x-x_{0}|\leq d(t)$ and $\rho(t,x_{0})=\underline{\rho}$. Here we choose $\beta=\frac{\gamma+1}{2}>0$,
\begin{equation}
\begin{aligned}
&{\rm supp}_{x\in{\rm supp}(\omega)}\e\rho^{\beta}(t,x)\leq\e\underline{\rho}^{\beta}+{\rm supp}_{x\in{\rm supp}(\omega)\cap A}\e\rho^{\beta}(t,x),\\
&\leq 2\e\underline{\rho}^{\beta}+{\rm supp}_{x\in{\rm supp}(\omega)\cap A}|\e\rho^{\beta}(t,x)-\e\rho^{\beta}(t,x_{0})|,\\
&\leq 2\e\underline{\rho}^{\beta}+{\rm supp}_{x\in{\rm supp}(\omega)\cap A}\int^{x_{0}+d(t)}_{x_{0}-d(t)}|\beta||\e\rho^{\beta-1}(t,x)\rho_{x}|dx,\\
&\leq 2\e\underline{\rho}^{\beta}+\int^{b+1+2d(t)}_{a-1-2d(t)}|\beta|\rho^{2\beta-1}dx+\int_{\R}\e^{2}\rho^{-1}\rho_{x}^{2}dx,\\
&\leq C(t)+\int^{b+1+2d(t)}_{a-1-2d(t)}\rho^{\gamma}dx,\\
&\leq C(t).
\end{aligned}
\label{2.28}
\end{equation}
Using (\ref{2.28}), Young inequalities and H\"older's inequalities, the first term on the right hand side of (\ref{2.25}) is treated as follows:
\begin{equation}
\begin{aligned}
&\int^{t}_{0}\int_{\R}\rho^{2}v_{x}\omega^{2}dx d\tau\\
&\leq\frac{1}{2}\e\int^{t}_{0}\int_{\R}\rho^{3}\omega^{4}dx d\tau+\frac{1}{2}\e\int^{t}_{0}\int_{\R}\rho v_{x}^{2}dx d\tau,\\
&\leq C(t)+\e\int^{t}_{0}\int_{\R}\rho^{3}\omega^{2}dx d\tau,\\
&\leq C(t)+C(t)\int^{t}_{0}\int_{\R}\rho^{4-\beta}\omega^{2}dx d\tau,\\
&\leq C(t)+\delta\int^{t}_{0}\int_{\R}\rho^{\gamma+1}\omega^{2}dx d\tau,\\
\end{aligned}
\label{2.29}
\end{equation}
Here we have used the fact that $\gamma>\frac{5}{3}$. By lemma \ref{lemma1} and the H\"older inequality, we obtain
\begin{equation}
\begin{aligned}
\big|\int^{x}_{-\infty}\rho v\omega dy\big|&\leq\int_{{\rm supp}(\omega)}|\rho v| dy,\\
&\leq (\int_{{\rm supp}(\omega)}\rho  dy)^{\frac{1}{2}}(\int_{{\rm supp}(\omega)}\rho v^{2} dy)^{\frac{1}{2}}\leq C(t).
\end{aligned}
\label{2.30}
\end{equation}
Then:
\begin{equation}
\begin{aligned}
&\big|\int_{\R}\big(\rho\omega\int^{x}_{-\infty}\rho v\omega  dy\big)dx\big|+\big|\int_{\R}\big(\rho_{0}\omega\int^{x}_{-\infty}\rho_{0} v_{0}\omega  dy\big)dx\big|\\
&\hspace{4cm}+\big|\int_{0}^{t}\int_{\R}\big(\rho u\omega_{x}\int^{x}_{-\infty}\rho v\omega  dy\big)dx d\tau\big|\leq C(t).
\end{aligned}
\label{2.32}
\end{equation}
Similarly, we have:
\begin{equation}
\big|\int_{0}^{t}\int_{\R}\big(\rho \omega\int^{x}_{-\infty}(\rho u v+P(\rho))\omega_{x}  dy\big)dx d\tau\big|\leq C(t),
\label{2.33}
\end{equation}
and
\begin{equation}
\begin{aligned}
&\e\big|\int_{0}^{t}\int_{\R}\big(\rho \omega\int^{x}_{-\infty}\rho v_{x}\omega_{x}  dy\big)dx d\tau\big|\\
&\leq \e\big|\int_{0}^{t}\int_{\R}\big(\rho \omega\int_{\R}\rho |v_{x}|\,|\omega_{x}|  dy\big)dx d\tau\big|,\\
&\leq \e\big|\int_{0}^{t}\big(\int_{\R}\rho \omega dx\big)\big(\int_{\R}\rho v_{x}^{2}dy+\int_{\R}\rho\omega_{x}^{2}dy\big)d\tau\big|,\\
&\leq C(t).
\end{aligned}
\label{2.34}
\end{equation}
Substituting (\ref{2.29}), (\ref{2.32})-(\ref{2.34}) into (\ref{2.25}) and noticing the smallness of $\delta$, we proved lemma \ref{lemme2}.
\begin{lem}
Suppose that  $(\rho_{0}(x),v_{0}(x)$ satisfy the conditions in the lemmas \ref{lemma1}. Furthermore there exists $M_{0}>0$ independent of $\e$, such that
\begin{equation}
\int_{\R}\rho_{0}(x)|v_{0}(x)-\bar{v}(x)|dx\leq M_{0}<+\infty,
\label{2.35}
\end{equation}
then for any compact set $K\subset\R$, it holds that:
\begin{equation}
\int^{t}_{0}\int_{K}(\rho^{\gamma+\theta}+\rho|v|^{3})dxd\tau\leq C(t,K),
\label{2.36}
\end{equation}
where $C(t,K)$ is independent of $\e$.
\end{lem}
\begin{remarka}
In order to prove the inequality (\ref{2.36}), we will use the same ingredients than in \cite{35} where this inequality was obtained for the first time.
\end{remarka}
{\bf Proof.}
We are now working with the function $\psi$ of proposition \ref{pair35}. If we consider $\eta^{\psi}_{m}$  as a function depending of $(\rho,v)$, we have for all $\rho\geq 0$ and $v\in\R$:
\begin{equation}
\begin{cases}
\begin{aligned}
&|\eta^{\psi}_{mv}(\rho,v)|\leq C,\\
&|\eta^{\psi}_{m\rho}(\rho,v)|\leq C\rho^{\theta-1}.
\end{aligned}
\end{cases}
\label{2.38}
\end{equation}
For this weak entropy pair $(\eta^{\psi},H^{\psi})$, we observe that:
$$\eta^{\psi}(\rho,0)=\eta^{\psi}_{\rho}(\rho,0)=0,\;H^{\psi}(\rho,0)=\frac{\theta}{2}\rho^{3\theta+1}\int_{\R}|s|^{3}[1-s^{2}]^{\lambda}_{+},$$
and:
$$\eta^{\psi}_{m}(\rho,0)=\beta\rho^{\theta}\;\;\mbox{with}\;\beta=\int_{\R}|s|[1-s^{2}]^{\lambda}_{+}ds.$$
By Taylor formula, we have:
\begin{equation}
\eta^{*}(\rho,m)=\beta\rho^{\theta}m+r(\rho,m),
\label{2.39}
\end{equation}
with:
\begin{equation}
r(\rho,m)\leq C\rho v^{2},
\label{2.40}
\end{equation}
for some constant $C>0$. Now we introduce a new entropy pair $(\hat{\eta},\hat{H})$ such that,
$$\hat{\eta}(\rho,m)=\eta^{\psi}(\rho,m-\rho v^{-}),\;\;\hat{H}(\rho,m)=H^{\psi}(\rho,m-\rho v^{-})+v^{-}\eta^{\psi}(\rho,m-\rho v^{-}),$$
with $m=\rho v$ which satisfies:
\begin{equation}
\begin{cases}
\begin{aligned}
&\hat{\eta}(\rho,m)=\beta\rho^{\theta+1}(v-v^{-})+r(\rho,\rho(v-v^{-})),\\
&r(\rho,\rho(v-v^{-}))\leq C \rho(v-v^{-})^{2}.
\end{aligned}
\end{cases}
\label{2.41}
\end{equation}
Integrating $(\ref{1.1})_{1}\times\hat{\eta}_{\rho}+(\ref{1.1})_{2}\times\hat{\eta}_{m}$ over $(0,t)\times(-\infty,x)$, we have:
\begin{equation}
\begin{aligned}
&\int^{x}_{-\infty}\big(\hat{\eta}(\rho,m)-\hat{\eta}(\rho_{0},m_{0})\big)dy+\int^{t}_{0}q^{*}(\rho,\rho(v-v^{-}))+v^{-}\eta^{*}(\rho,\rho(v-v^{-}))d\tau\\
&=t q^{*}(\rho^{-},0)+\e \int^{t}_{0}\hat{\eta}_{m}\rho  v_{x} d\tau-\e\int^{t}_{0}\int^{x}_{-\infty}(\hat{\eta}_{mu}\rho v_{x}^{2}
+\hat{\eta}_{m\rho}\rho \rho_{x}v_{x})dyd\tau.
\end{aligned}
\label{2.42}
\end{equation}
By using (\ref{2.38}), we obtain:
\begin{equation}
\big|\e\int^{t}_{0}\int^{x}_{-\infty}\hat{\eta}_{mu}\rho v_{x}^{2}dyd\tau\big|\leq C\e\int^{t}_{0}\int_{\R}\rho v_{x}^{2}dy\,d\tau\leq C(t),
\label{2.43}
\end{equation}
\begin{equation}
\begin{aligned}
&\big|\e\int^{t}_{0}\int^{x}_{-\infty}\hat{\eta}_{m\rho}\rho \rho_{x}v_{x}dyd\tau\big|\leq C\e\int^{t}_{0}\int_{\R}\rho^{\theta-1}\rho  |\rho_{x} v_{x}|dy\,d\tau\leq C(t),\\
&\hspace{2cm}\leq C\e\int^{t}_{0}\int_{\R}\rho v_{x}^{2}dy\,d\tau+C\e\int^{t}_{0}\int_{\R}\rho^{\gamma-2}\rho_{x}^{2}dy\,d\tau\leq C(t).
\end{aligned}
\label{2.44}
\end{equation}
Substituting (\ref{2.43}) and (\ref{2.44}) into (\ref{2.42}), then integrating over $K$ and using (\ref{2.37}), we obtain:
\begin{equation}
\begin{aligned}
&\int^{t}_{0}\int_{K}\rho^{\theta+\gamma}+\rho|v-v^{-}|^{3}dxd\tau\\
&\leq C(t)+C\int^{t}_{0}\int_{K}|\eta^{*}(\rho,\rho(v-v^{-})|dxd\tau+C\e
\int^{t}_{0}\int_{K}\rho|v||v_{x}|dxd\tau\\
&+C\e \int^{t}_{0}\int_{K}\rho^{1+\theta}|v_{x}|dxd\tau+2\sup_{\tau\in [0,t]}\big|\int_{K}(\int^{x}_{-\infty}\hat{eta}(\rho(y,\tau),(\rho v)(y,\tau))dy)dx\big|.
\end{aligned}
\label{2.45}
\end{equation}
Applying lemma \ref{lemma1}, we have:
\begin{equation}
\int^{t}_{0}\int_{K}|\eta^{*}(\rho,\rho(v-v^{-})|dxd\tau\leq C(t).
\label{2.46}
\end{equation}
By H\"older's inequality and (\ref{2.28}), we get:
\begin{equation}
\begin{aligned}
\e \int^{t}_{0}\int_{K}\rho^{1+\theta}|v_{x}|dxd\tau&\leq C\e\int^{t}_{0}\int_{K}\rho v_{x}^{2}dxd\tau
+C\e\int^{t}_{0}\int_{K}\rho^{1+2\theta}dxd\tau,\\
&\leq C(t)+C(t)\int^{t}_{0}\int_{K}\rho^{\theta}dxd\tau,\\
&\leq C(t).
\end{aligned}
\label{2.47}
\end{equation}
We have now:
\begin{equation}
\begin{aligned}
\e \int^{t}_{0}\int_{K}\rho|v||v_{x}|dxd\tau&\leq\frac{1}{2}\e\int^{t}_{0}\int_{K}\rho v_{x}^{2}dx d\tau+
\frac{1}{2}\e\int^{t}_{0}\int_{K}\rho v^{2}dx d\tau,\\
&\leq C(t).
\end{aligned}
\label{2.48}
\end{equation}
Now we are going to deal with the last term on the right hand side of (\ref{2.45}). (\ref{1.1}) implies that:
\begin{equation}
(\rho v-\rho v^{-})_{t}+(\rho v^{2}+P(\rho)-\rho u  u^{-})_{x}=\e(\rho v_{x})_{x}.
\label{2.49}
\end{equation}
Integrating (\ref{2.49}) over $[0,t]\times (-\infty,x)$ for $x\in K$, we get:
\begin{equation}
\begin{aligned}
&\int^{x}_{-\infty}\rho (v-\rho v^{-})dy=\int^{x}_{-\infty}\rho_{0} (v_{0}-\rho v^{-})dy-\int^{t}_{0}(\rho v^{2}+P(\rho)-\rho u  u^{-}-P(\rho^{-}))\\
&\hspace{10cm}+\e\int^{t}_{0}\rho v_{x}d\tau.
\end{aligned}
\label{2.50}
\end{equation}
Furthermore:
\begin{equation}
\begin{aligned}
&\big|\int^{x}_{-\infty}\hat{\eta}((\rho(y,\tau),(\rho v)(y,\tau))dy\big|\\
&\leq\big|\int^{x}_{-\infty}(\hat{\eta}(\rho\rho v)-\beta\rho^{\theta+1}(v-\bar{v}))dy\big|+
|\int^{x}_{-\infty}\beta\rho^{\theta+1}(v-\bar{v}))dy\big|\\
&\leq\big|\int^{x}_{-\infty}(r(\rho\rho (v-\bar{v}))dy\big|+\big|\int^{x}_{-\infty}\beta(\rho^{\theta}-(\rho^{-})^{\theta})\rho(v-\bar{v}))dy\big|\\
&\hspace{8cm}+
\beta(\rho^{-})^{\theta}\big|\int^{x}_{-\infty}\rho(v-\bar{v}))dy\big|,\\
&\leq C(t)+
\beta(\rho^{-})^{\theta}\big|\int^{x}_{-\infty}\rho(v-\bar{v}))dy\big|.
\end{aligned}
\label{2.51}
\end{equation}
By using (\ref{2.35}), lemma \ref{lemma1} and \ref{lemme2}, (\ref{2.50}) and (\ref{2.51}) we conclude the proof of the lemma.
\subsection{$H^{-1}_{loc}(\R^{2}_{+})$ Compactness}
In this section we are going to take profit of the uniform estimates obtained in the previous section in order to prove the following lemma, which gives the $H^{-1}_{loc}(\R^{2}_{+})$-compactness of the Korteweg solution sequence $(\rho^{\e},v^{\e})$ on a entropy- entropy flux pair.
\begin{lem}
Let $\psi\in C^{2}_{0}(\R)$, $\eta^{\psi},H^{\psi})$ be a weak entropy pair generated by $\psi$. Then for the solutions $(\rho^{\e},v^{\e})$ with $m^{\e}=\rho^{\e}v^{\e}$ of Korteweg system (\ref{1.1}, the following sequence:
\begin{equation}
\eta^{\psi}(\rho^{\e},m^{\e})_{t}+q^{\psi}(\rho^{\e},m^{\e})_{x}\;\;\mbox{are compact in}\;H^{-1}_{loc}(\R^{2}_{+})
\label{3.1}
\end{equation}
\label{lemme4}
\end{lem}
{\bf Proof:} Now we are going to prove the lemma. A direct computation on $(\ref{1.1})_{1}\times\eta_{\rho}^{\psi}(\rho^{\e}, m^{\e})+(\ref{1.1})_{2}\times\eta_{m}^{\psi}(\rho^{\e}, m^{\e})$ gives:
\begin{equation}
\begin{aligned}
&\eta_{\rho}^{\psi}(\rho^{\e}, m^{\e})_{t}+H_{\rho}^{\psi}(\rho^{\e}, m^{\e})_{x}=\e\big(\eta_{\rho}^{\psi}(\rho^{\e}, m^{\e})(\rho^{\e}) v^{\e}_{x}\big)-\e\eta^{\psi}_{mu}(\rho^{\e}, m^{\e})(\rho^{\e})(v^{\e}_{x})^{2}\\
&\hspace{8cm}-\e\eta^{\psi}_{mu}(\rho^{\e}, m^{\e})(\rho^{\e}) v^{\e}_{x}\rho^{\e}_{x}.
\end{aligned}
\label{3.7}
\end{equation}
Let $K\subset\R$ be compact, using proposition \ref{propChen} (\ref{3.6}) and H\"older inequality, we get:
\begin{equation}
\begin{aligned}
&\e\int^{t}_{0}\int_{K}|\eta^{\psi}_{mu}(\rho^{\e}, m^{\e})(\rho^{\e})|(v^{\e}_{x})^{2}+
|\eta^{\psi}_{mu}(\rho^{\e}, m^{\e})(\rho^{\e})v^{\e}_{x}\rho^{\e}_{x}|dxdt \\
&\leq C\e\int^{t}_{0}\int_{K}(\rho^{\e})|(v^{\e}_{x})^{2}dxd\tau+
C\e\int^{t}_{0}\int_{K}(\rho^{\e})^{\gamma-2}(\rho^{\e}_{x})^{2}dxd\tau\\
&\leq C(t).
\end{aligned}
\label{3.8}
\end{equation}
This shows that:
\begin{equation}
-\e\eta^{\psi}_{mu}(\rho^{\e}, m^{\e})\rho^{\e}(v^{\e}_{x})^{2}-\e\eta^{\psi}_{mu}(\rho^{\e}, m^{\e})\rho^{\e} v^{\e}_{x}\rho^{\e}_{x}\;\;\mbox{are bounded in}\; L^{1}([0,T]\times K),
\label{3.9}
\end{equation}
and thus it is compact in $W^{-1,p_{1}}_{loc}(\R^{2}_{+})$, for $1<p_{1}<2$. Moreover we observe that
$$|\eta^{\psi}_{mu}(\rho^{\e}, \rho^{\e}v^{\e})|\leq C_{\psi},$$,
then we obtain:
\begin{equation}
\begin{aligned}
&\int^{t}_{0}\int_{K}(\e\eta^{\psi}_{m}(\rho^{\e},m^{\e})\rho^{\e}v_{x}^{\e})^{\frac{4}{3}}dxdt\\
&\leq \int^{t}_{0}\int_{K}\e^{\frac{4}{3}}|\rho^{\e}|^{\frac{4}{3}}|v_{x}^{\e}|^{\frac{4}{3}}dxdt\\
&\leq C\e^{\frac{4}{3}} \int^{t}_{0}\int_{K}\rho^{\e}|v_{x}^{\e}|^{2}dxdt+
 C\e^{\frac{4}{3}} \int^{t}_{0}\int_{K}(\rho^{\e})^{2}dxdt\\
&\leq C(t,K) \e^{\frac{1}{3}}+ C\e^{\frac{4}{3}} \int^{t}_{0}\int_{K}(\rho^{\e})^{\gamma+1}dxdt\rightarrow_{\e\rightarrow 0}0.
\end{aligned}
\label{3.10}
\end{equation}
Using (\ref{3.10}) and (\ref{3.9}), we obtain
\begin{equation}
\eta_{\rho}^{\psi}(\rho^{\e}, m^{\e})_{t}+H_{\rho}^{\psi}(\rho^{\e}, m^{\e})_{x}\;\mbox{are compact in}\; W^{-1,p_{1}}_{loc}(\R^{2}_{+})\;\mbox{for some}\;1<p_{1}<2.
\label{3.11}
\end{equation}
Furthermore by (\ref{3.3})-(\ref{3.4}), lemma \ref{lemma1}-\ref{lemme2} and (\ref{2.36}), we have:
\begin{equation}
\eta_{\rho}^{\psi}(\rho^{\e}, m^{\e})_{t}+H_{\rho}^{\psi}(\rho^{\e}, m^{\e})_{x}\;\mbox{are uniformly bounded in}\; L^{p_{3}}_{loc}(\R^{2}_{+})\;\mbox{for}\;p_{3}>2.
\label{3.12}
\end{equation}
where $p_{3}=\gamma+1>2$ when $\gamma\in(1,3]$, and $p_{3}=\frac{\gamma+\theta}{1+\theta}>2$ when $\gamma>3$. By interpolation we conclude the proof of the lemma \ref{lemme4}.
\section{Proof of theorem \ref{theo1}}
From lemmas \ref{lemma1}, we have verified the conditions (i)-(iii) of theorem \ref{theo2} for the sequence of solutions $(\rho^{\e},m^{\e})$. 
Using theorem \ref{theo2}, there exists a subsequence $(\rho^{\e},m^{\e})$ and a pair of measurable functions $(\rho,m)$ such that
\begin{equation}
(\rho^{\e},m^{\e})\rightarrow (\rho,m),\;\;\mbox{a.e}\;\e\rightarrow 0.
\label{4.1}
\end{equation}
It is easy to check that $(\rho,m)$ is a finite-energy entropy solution $(\rho,m)$ to the Cauchy problem (\ref{1.3}) with initial data $(\rho_{0},\rho_{0}u_{0})$ for the isentropic Euler equations with $\gamma>\frac{5}{3}$. It achieves the proof of theorem \ref{theo1}.


\begin{thebibliography}{00}
\bibitem{fA}
D.M. Anderson, G.B McFadden and A.A. Wheller. Diffuse-interface
methods in
fluid mech. \textit{In Annal review of fluid mechanics}, Vol. 30, pages 139-165. Annual Reviews, Palo Alto,
CA, 1998.
\bibitem{BB1}
S. Bianchini and A. Bressan, Vanishing viscosity solutions of nonlinear hyperbolic systems. \textit{Ann. Math.},  {\bf 161} 1 (2005), pp. 223Ð342.
\bibitem{fC}
J.W. Cahn, J.E. Hilliard, Free energy of a nonuniform system, I.
Interfacial free energy, \textit{J. Chem. Phys.} 28 (1998) 258-267.
\bibitem{Chen}
G. Q. Chen, The theory of compensated compactness and the system of isentropic gas dynamics, \textit{Lecture notes, Preprint MSRI-00527-91}, Berkeley, October 1990.
\bibitem{10}
G. Chen and M. Perepelista, Vanishing viscosity limit of the Navier-Stokes equations to the Euler equations for compressible fluid flow, \textit{Comm. Pure. Appl. Math.} 2010 (to appear).
\bibitem{fDD}
R. Danchin and B. Desjardins, Existence of solutions for
compressible fluid models of Korteweg type, \textit{Annales de l'IHP, Analyse
non
lin\'eaire} 18,97-133 (2001).
\bibitem{Di1}
R. J. DiPerna, Convergence of the viscosity method for isentropic gas dynamics, \textit{Commun. Math. Phys.} 91 (1983), 1-30.
\bibitem{Di2}
R. J. DiPerna, Convergence of approximate solutions to conservation laws, \textit{Arch. Rat. Mech. Anal.}  (1983), 22-70.
\bibitem{fDS}
J.E. Dunn and J. Serrin, On the thermomechanics of interstitial
working , \textit{Arch. Rational Mech. Anal}. 88(2) (1985) 95-133.
\bibitem{Evans}
L. C. Evans, Partial differential equations, \textit{AMS}, GSM Volume 19.
\bibitem{G}
J. Glimm, Solutions in the large for nonlinear hyperbolic systems of equations,.\textit{Comm. Pure Appl. Math. 18}. 1965 (697-715).
\bibitem{fGP}
M.E. Gurtin, D. Poligone and J. Vinals, Two-phases binary fluids and
immiscible fluids described by an order parameter, \textit{Math. Models
Methods Appl. Sci}. 6(6) (1996) 815--831.
\bibitem{Hprepa}
B. Haspot, Blow-up criterion, ill-posedness and existence of  strong solution for Korteweg system with infinite energy, \textit{preprint}  (arxiv February 2011).
\bibitem{Hprepa1}
B. Haspot, New entropy for KortewegÕs system, existence of global weak solution and Prodi-Serrin theorem, \textit{preprint}  (arxiv February 2011).
\bibitem{fH1}
B. Haspot, Existence of solutions for compressible fluid models of Korteweg type, \textit{Annales Math\'ematiques Blaise Pascal} {\bf 16}, 431-481 (2009).
\bibitem{fH2}
B. Haspot, Existence of weak solution for compressible fluid models of Korteweg type, \textit{Journal of Mathematical Fluid Mechanics}, DOI: 10.1007/s00021-009-0013-2 online.
\bibitem{HM}
M. Heida and J. M\'alek, On compressible Korteweg fluid-like materials, \textit{International Journal of Engineering Science}, Volume 48, Issue 11, November 2010, Pages 1313-1324.
\bibitem{Hu1}
F. Huang and Z. Wang, Convergence of viscosity solutions for isothermal gas dynamics, \textit{SIAM J. Math. Anal.} 34 (2002), no3, 595-610.
\bibitem{Hu2}
F. Huang, R. Pan, T. Wang, Y. Wang and and X. Zhai, Vanishing viscosity limit for isentropic Navier-Stokes equations with density-dependent viscosity. \textit{Preprint}.
\bibitem{fJL}
D. Jamet, O. Lebaigue, N. Coutris and J.M. Delhaye. The second
gradient method for the direct numerical simulation of liquid-vapor
flows
with phase change. \textit{J. Comput. Phys}, 169(2): 624--651, (2001).
\bibitem{Lax}
P. Lax, Hyperbolic systems of conservation laws II. \textit{Comm. on Pure and Applied Math.}, 10: 537-566, 1957.
\bibitem{fK}
D.J. Korteweg. Sur la forme que prennent les \'equations du
mouvement des fluides si l'on tient compte des forces capillaires
par des variations de densit\'e. \textit{Arch. N\'eer. Sci. Exactes S\'er}.
II, 6 :1-24, 1901.
\bibitem{Lefloch}
 P.G. LeFloch, Hyperbolic Systems of Conservation Laws. The theory of classical and nonclassical shock waves, \textit{Lectures in Mathematics}, ETH Z\"urich, Birkh\"auser, 2002.
\bibitem{fL2}
P.-L. Lions, Mathematical Topics in Fluid Mechanics, Vol 2,
Compressible models, \textit{Oxford University Press} (1996).
\bibitem{36}
P.-L. Lions, B. Perthame and P. Souganidis, Existence and stability of entropy solutions for the hyperbolic systems of isentropic gas dynamics in Eulerian and Lagrangian coordinates, \textit{Comm. Pure. Apll. Math.} 49 (1996), 599-638.
\bibitem{35}
P.-L. Lions, B. Perthame and E. Tadmor Mathematical, Kinetic formulation of the isentropic gas dynamics and p-systems, \textit{Commun. Phy.} 163 (1994), 415-431.
\bibitem{MV}
A.Mellet and A.Vasseur, Existence and uniqueness of global strong solutions for one-dimensional compressible Navier-Stokes equations, \textit{SIAM J. Math. Anal. 39} (2008), No 4, 1344-1365.
\bibitem{REF}
L. Modica, The gradient theory of phase transitions and the minimal interface criterion. \textit{Arch. Rational Mech. Anal}. {\bf 98}, no. 2, 123-142 (1987).
\bibitem{Mu}
F. Murat, Compacit\'e par compensation, \textit{Ann. Sci. Norm. Sup. Pisa. 5}, 1978, p 489-507.
\bibitem{Per}
B. Perthame, Higher moments lemma: application to Vlasov-Poisson and Fokker-Planck equations. \textit{Math. Meth. Appl. Sc.} 13, 441-452 (1990).
\bibitem{Rohdehdr}
C. Rohde,   Approximation of Solutions of Conservation Laws by Non- Local Regularization and Discretization, \textit{Habilitation Thesis}, University of Freiburg (2004).
\bibitem{fR}
J.S. Rowlinson, Translation of J.D van der Waals. The thermodynamic
theory of capillarity under the hypothesis of a continuous variation
of density. \textit{J.Statist. Phys.}, 20(2): 197-244, 1979.
\bibitem{Serre}
D. Serre, Systems of conservation laws, I- Hyperbolicity, Entropies, Shock Waves, \textit{ Cambridge University Press} (1999).
\bibitem{Sol}
V. A Solonnikov. The solvability of the initial boundary-value problem for the equations of motion of a viscous compressible fluid. \textit{Zap Naucn. Sem. Leningrad. Odel. Mat. Inst. Steklov.} (LOMI), 56: 128-142, 197, 1976. Investigations on linear operators and theory of functions, VI.
\bibitem{Ta}
L. Tartar, Compensated compactness and applications to partial differential equations,  \textit{Nonlinear Analysis and Mechanics , Herriot-Watt Symposium}, Vol 4, R. J. Knops, ed, Research Notes in Mathematics No 39, Pitman, Boston-London, 1979, p 136-212.
\bibitem{Ta1}
L. Tartar,  The compensated compactness method applied to systems of conservation laws, \textit{Systems of Nonlinear Partial Differential Equations, NATO Adv. Sci. Inst. Ser. C: Math. Phys. Sci. 111}, J. Bali, ed, Reidel, Dordrecht-Boston, 1983, p 263-285.
\bibitem{fTN}
C. Truedelland W. Noll. The nonlinear field theories of mechanics.
\textit{Springer-Verlag}, Berlin, second edition, 1992.
\bibitem{VW}
J.F Van der Waals. Thermodynamische Theorie der Kapillarit\"at unter Voraussetzung stetiger Dichte\"anderung, \textit{Phys. Chem.} {\bf 13}, 657-725 (1894).
\end{thebibliography}
\end{document}